\RequirePackage[l2tabu, orthodox]{nag}
\RequirePackage{snapshot}
\documentclass[journal, 10pt]{amsart}

\usepackage{srcltx}
\usepackage{array}

\usepackage{graphicx}
\usepackage{subfigure} 

\usepackage[utf8x]{inputenc}
\usepackage[english]{babel}

\usepackage{amssymb, amsmath}
\usepackage{bm}

\usepackage{epsfig}
\usepackage{color}
\usepackage{psfrag}

\usepackage{fullpage}
\usepackage{setspace}
\usepackage{flushend}

\usepackage{algorithmic}
\usepackage[ruled]{algorithm2e}

\usepackage[short,12hr]{datetime}

\usepackage{booktabs}
\usepackage{multirow}

\usepackage{cuted}
\usepackage{latexsym}

\usepackage{accents}

\newtheorem{theorem}{Theorem}

\newtheorem{corollary}{Corollary}

\makeatletter

\newenvironment{rsmallmatrix}{\null\,\vcenter\bgroup
  \Let@\restore@math@cr\default@tag
  \baselineskip6\ex@ \lineskip1.5\ex@ \lineskiplimit\lineskip
  \ialign\bgroup\hfil$\m@th\scriptstyle##$&&\thickspace\hfil
  $\m@th\scriptstyle##$\crcr
}{%
  \crcr\egroup\egroup\,%
}
\makeatother

\DeclareSymbolFont{epsilon}{OML}{ntxmi}{m}{it}
\DeclareMathSymbol{\epsilon}{\mathord}{epsilon}{"0F}

\title{On the Computation of Neumann Series}

\author{%
Vassil~S.~Dimitrov
\qquad
Diego~F.~G.~Coelho}

\thanks{V.~S.~Dimitrov is with
the Department of Electrical and Computer Engineering, University of Calgary, Calgary, AB, Canada
and
the Computer Modelling Group, Ltd., AB, Canada.
}%

\thanks{
D.~F.~G.~Coelho
is with
the Department of Electrical and Computer Engineering, University of Calgary, Calgary, AB, Canada.
}

\thanks{This project was supported by the  Natural Sciences and Engineering Research Council of Canada (NSERC)}

\begin{document}

\begin{abstract}
This paper proposes new factorizations for computing the Neumann series.
The factorizations are based on fast algorithms for small prime sizes series and the splitting of large sizes into several smaller ones.
We propose a different basis for factorizations other than the well-known binary and ternary basis.
We show that is possible to reduce the overall complexity for the usual binary decomposition from~$2\operatorname{log}_2(N)-2$ multiplications to around~$1.72\operatorname{log}_2(N)-2$ using a basis of size five.
Merging different basis we can demonstrate that we can build fast algorithms for particular sizes.
We also show the asymptotic case where one can reduce the number of multiplications to around~$1.70\operatorname{log}_2(N)-2$.
Simulations are performed for applications in the context of wireless communications and image rendering, where is necessary perform large sized matrices inversion.
\end{abstract}


\maketitle

\section{Introduction}

Neumann series was introduced by Carl Neumann in 1877 in connection with function analysis.
%
%
Neumann series has been applied to theoretical problems for solving Fredholm integral equation.
However, apart from theoretical application of Neumann series, it has been overlooked by scientific community in the application for solving computational problems.

Neumann series has been used in contexts other than function analysis where the computation of Neumann series for certain matrices are necessary~\cite{Dimitrov1995a}.
It has been applied to image rendering and wireless communication~\cite{Ng2011, Prabhu2014} for computing  approximate inverse matrices.
In both scenarios, the Neumann series computation has been implemented straightforwardly by applying Horner rule~\cite{Graham1989}.
Neumann series also occur in the inversion of elements on finite fields~$GF(2^m)$~\cite{Jarvinen2014,Azarderakhsh2012,Liu2014} and quantum computing~\cite{Kalapala2008}.
Despite these applications, little effort has been done to develop efficient algorithms for computing the Neumann series of matrices.
Usual approaches for speedup the Neumann series computation for matrices has been based on binary decomposition of the associated polynomial used for approximating the Neumann series~\cite{Graham1989}.
To the best of our knowledge, just a few papers has been devoted to develop algorithms for Neumann series evaluation.

The first approach occurred in literature for minimizing the number of operations for evaluating Neumann series is due to Westreich~\cite{Westreich1989}.
In~\cite{Westreich1989}, Westreich proposes a method for minimizing the number of matrix additions and multiplications achieving~$3\operatorname{log}_2(N)$ matrix multiplications.
Roy and Minocha proposed a scheme that also requires around~$3\operatorname{log}_2(N)$ matrix multiplications in~\cite{Roy1992}.
Lei and Nakamura proposed the well-known binary approach for Neumann series~\cite{Lei1992} for reducing matrix multiplications and thus achieving~$2\operatorname{log}_2(N)-2$ matrix multiplications.
Lei and Nakamura also conjectured that no better algorithm could be found.
However, in~\cite{Dimitrov1994}, Dimitrov and Donevsky  provided a disproof of such a conjecture by providing a factorization for different Neumann series sizes requiring less matrix multiplications than~$2\operatorname{log}_2(N)-2$.
In~\cite{Dimitrov1995a}, Dimitrov and Cooklev have proposed a hybrid algorithm considering ternary decomposition and binary decomposition, requiring around~$1.94\operatorname{log}_2(N)-2$ matrix multiplications
%

In this paper we develop efficient algorithms for sizes that are power of a small prime number other than two or three.
We show that the previous methods are not optimal and we can reduce the usual binary approach complexity from~$2\operatorname{log}_2(N)-2$ matrix multiplications to about~$1.72\operatorname{log}_2(N)-2$ using basis five.
We provide a general theory for mixing approaches based on different basis and we provide asymptotic case multiplicative complexity achieving around~$1.701\operatorname{log}_2(N)-2$.
Experiments in the context of image rendering and wireless communication~\cite{Ng2011, Prabhu2014} are provided along with simulation results.

The paper unfolds as the following.
Section~\ref{sec:math} provides the background for the theory developed on the paper.
Section~\ref{sec:decomp-series-prime} shows the decomposition of big length sequences that are powers of prime numbers.
Algorithms for small prime lengths are provided along with the complexity analysis.
Section~\ref{sec:optimal-base} analyzes the asymptotic case where a optimal basis is found using recurrence series.
Complexity analysis are shown using results provided by Aho and Sloan~\cite{Aho1973}.
In Section~\ref{sec:base-combination}, we mix different basis in order to take advantage of each basis for the computation of Neumann series.
In Section~\ref{sec:app} we apply the theory developed for reducing matrix multiplication in the context of image rendering and wireless communication.
We run simulations with matrices of several sizes and we approximate their inverse.
Conclusions and future works comments are drawn on Section~\ref{sec:conclusion}.

\section{Mathematical Background}
\label{sec:math}

The Neumann series of a matrix is an infinity series holding the following equality
\begin{align*}
\mathbf{A}^{-1}  = \sum_{n = 0}^{\infty} (\mathbf{I}-\mathbf{A})^n,
\end{align*}
where the matrix~$\mathbf{A}$ has eigenvalues with real part on the interval~$(0,2)$.
The above series can be approximated with an arbitrary number of terms~$N$.
The number of terms used for approximating the inverse of~$\mathbf{A}$ depends on the context application and naturally on how well conditioned is the matrix.
A numerical evaluation of~$\mathbf{A}^{-1}$ can be seen as the computation of a geometric series with a finite number of terms as
\begin{align}
\label{eq:invAN}
\mathbf{A}^{-1}_N  \approx \sum_{n = 0}^{N-1} \mathbf{B}^n,
\end{align}
where~$\mathbf{B} = \mathbf{I}-\mathbf{A}$ and~$\mathbf{A}^{-1}_N$ represents the~$\mathbf{A}^{-1}$ approximation with only~$N$ terms.
%
%
Note that the right-hand side of the expression in~\eqref{eq:invAN} is equivalent to a geometric series with a real quantity~$x$.
We denote a geometric series on variable~$x$ of~$N$ terms as the following.
Let~$x$ be a real.
%
%
We define the geometric series on~$x$ with $N$ terms as the sum of consecutive powers of~$x$ as
$f(N, x) = \sum_{n = 0}^{N-1}x^n$.

\section{Geometric Series Decomposition into Prime Size Series}
\label{sec:decomp-series-prime}
Let~$N$ be a composite number.
A geometric series of size~$N$ can be decomposed as a product of several small geometric series.
This is guaranteed by the next theorem.
%
%
\begin{theorem}[Geometric Series Decomposition]
\label{theo:geometric-series-reduction}
Let~$f(N,x)$ be a geometric series with~$N$ terms in~$x$ and let~$N = KJ$. We have that
~$f(N,x) = f(K,x)\cdot f(J,x^K)$.

\begin{proof}
From the right side 
%
%
we have that
\begin{align*}
f(K,x)\cdot f(J,x^K) & = \left( \sum_{k = 0}^{K-1}x^k \right) \cdot \left( \sum_{j = 0}^{J-1} x^{K\cdot j} \right)\\
& = \sum_{j = 0}^{J-1} \left( x^{K\cdot j} \left( \sum_{k = 0}^{K-1}x^k \right)\right)\\
& = \sum_{j = 0}^{J-1} \left( \sum_{k = 0}^{K-1}x^{k+K\cdot j} \right).
\end{align*}
The inner summation varies on~$k$ and is nothing more than~$1+x+x^2+\ldots+x^{K-1}$.
For each step of outer summation that is in~$j$ we have a shift of~$K$ units in the exponent of each power of~$x$.
For~$j = 0$ we have the parcel~$1+x+x^2+\ldots+x^{K-1}$, for~$j=1$ we have parcel~$(1+x+x^2+\ldots+x^{K-1})\cdot x^K$ and so on until~$(1+x+x^2+\ldots+x^{K-1})\cdot x^{K(J-1)}$.
Thus, in the first parcel ($j=0$) we have power of~$x$ from~$0$ to~$K-1$, in the second parcel ($j=1$) we have powers of~$x$ from~$K$ to~$2K-1$ and so on until the~$K(J-1)$th parcel ($j=K(J-1)$), which results in powers of~$x$ from~$K(J-1)$ to~$KJ=N$.
Therefore, we have that
\begin{align*}
f(K,x)\cdot f(J,x^K) & = \sum_{j = 0}^{J-1} \left( \sum_{k = 0}^{K-1}x^{k+K\cdot j} \right)\\
& = \sum_{n = 0}^{N-1} x^n \\
& = f(N,x).
\end{align*}
\end{proof}
\end{theorem}

%
All the following development on this paper relies on Theorem~\ref{theo:geometric-series-reduction}.
All early works on the computation of Neumann series have concentrated on binary and ternary factorization~\cite{Stewart1998,Dimitrov1995a}.
%
%
In the following we generalize this for arbitrary power of prime length.

\subsection{Power of Primes Size}

Consider now the case where~$N$ is a power of a prime integer~$P$.
We use Theorem~\ref{theo:geometric-series-reduction} for the following corollary.
\begin{corollary}[Reduction for Geometric Series with Prime Terms]
\label{cor:reduction-geometric-series-prime}
Let~$N$ be a power of an integer~$P$.
The geometric series of~$N$ terms on~$x$,~$f(N,x)$, can be expressed as
\begin{equation*}
f(N, x) = \prod_{i=0}^{\operatorname{log}_{P}(N)-1} f(P, x^{P^i}),
\end{equation*}
where~$\operatorname{log}_{P}(N)$ is the logarithm of~$N$ in base~$P$.

\begin{proof}
Apply Theorem~\ref{theo:geometric-series-reduction} once to~$f(N, x)$, what leads to
\begin{equation*}
f(N, x) = f(P, x) \cdot f(P^{\operatorname{log}_{P}(N)-1}, x^P).
\end{equation*}
Makes the same for~$f(P^{\operatorname{log}_{P}(N)-1}, x^P)$, what leads to
\begin{equation*}
f(N, x) =  f(P, x) \cdot f(P, x^{P})\cdot f(P^{\operatorname{log}_{P}(N)-2}, x^{P^2}).
\end{equation*}
If we repeat the application of Theorem~\ref{theo:geometric-series-reduction} we obtain
\begin{align}
\label{eq:prime-case}
f(N, x) & =  f(P, x) \cdot f(P, x^{P}) \cdots f(P, x^{P^{\operatorname{log}_{P}(N)-1}})\nonumber\\
& = \prod_{i=0}^{\operatorname{log}_{P}(N)-1} f(P, x^{P^i}).
\end{align}
\end{proof}
\end{corollary}

Particular examples for~$P = 2, 3$ and~$5$, respectively, result in 
\begin{align*}
f(N, x) &  =  \prod_{i=0}^{\operatorname{log}_{2}(N)-1} \left( 1+x^{2^i} \right)\\
f(N, x) &  =  \prod_{i=0}^{\operatorname{log}_{3}(N)-1} \left( 1+x^{3^i} + x^{2\cdot 3^i} \right)\\
f(N, x) &  =  \prod_{i=0}^{\operatorname{log}_{5}(N)-1} \left( 1+x^{5^i} + x^{2\cdot 5^i} + x^{3\cdot 5^i} + x^{4\cdot 5^i} \right).
\end{align*}


In the following, we provide complexity analysis for the power of primes case when we employ direct evaluation of~$f(P, x^{P^i})$.
By direct evaluation we mean the use of Horner rule~\cite{Graham1989} for computing the Neumann series.
\begin{theorem}[Complexity of Direct Evaluation of Power of Prime Size Series]
\label{theo:complexity-direct-prime}
%
The direct evaluation of~$f(N, x)$ when~$N$ is a power of a prime~$P$ requires~$P\cdot \operatorname{log}_P(N)-2$ multiplications.

\begin{proof}
Consider the direct evaluation of~$f(N, x)$ in Corollary~\ref{cor:reduction-geometric-series-prime}.
%
%
In the first iteration, when~$i = 0$, we need~$P-2$ multiplications.
In order to compute~$x^P$ as required for the next iteration, we do~$x^P = f(P,x)\cdot (x-1)+1$, thus requiring one more multiplication.
In the second iteration, when~$i = 1$, we need more~$P-2$ multiplications plus one multiplication to compute~$x^{P^2}$ by doing~$x^{P^2} = f(P,x^P)\cdot (x^P-1)+1$.
For each iteration, we then need~$P-1$ multiplications.
In the last multiplication, we just need~$P-2$ multiplications, not~$P-1$.
This is because there is no next iteration.
Since we have~$\operatorname{log}_P(N)$ iterations, in order to evaluate all~$f(P, x^{P^i})$ parcels we need~$$(P-2) + (P-1)\cdot (\operatorname{log}_P(N)-1).$$
In order to calculate the final quantity~$f(N,x)$, we need more~$\operatorname{log}_P(N)-1$ to multiply the result in each iteration.
This results in a total of 
%
%
of~$P\cdot \operatorname{log}_P(N)-2$.
\end{proof}
\end{theorem}

%
If we evaluate the small prime series of size~$P$ directly, one could argue that Theorem~\ref{theo:complexity-direct-prime} shows that no one could find a smaller number of multiplications than the one offered by the binary and ternary factorization~\cite{Lei1992, Dimitrov1995a}.
In fact, if one take~$P$ other than two or three, we have that~$P\operatorname{log}_P(N) = P/\operatorname{log}_2(P) \operatorname{log}_2(N) > 2\operatorname{log}_2(N)$ always that~$P-2\operatorname{log}_2(P) > 0$.
This is true for all primes~$P \geq 5$.
This shows the importance of evaluation of geometric series with prime terms.
However, we will use Theorem~\ref{theo:complexity-direct-prime} in connection with smart factorization of prime length sequences to show that one can attain reduced complexity using power of primes other than two and three.
Therefore, we concentrate our following analysis in the evaluation of geometric series of prime size.

\subsection{Fast Algorithms for Prime Length Case}
\label{sec:fast-alg-small}
In the following, we develop fast algorithms for each small prime size and address the complexity for sequences that have size that is a power of a prime.
The evaluation for very small size is straightforward, such as for sizes~$P = 2$ or~$3$, which directly we compute~$f(2, x) = 1+x$ and~$f(3,x) = 1+x+x^2$.
But when the length of the sequence increases, smart factorization are necessary.

For size~$5$, for example, instead of compute~$f(5,x)$ directly, which would require~$3$ multiplication, one can note that~$f(5,x) = 1+(1+x^2)\cdot(x+x^2)$, which requires only~$2$ multiplications.
In order to evaluate~$x^5$ to be used for sequences of size power of five, one needs only one more multiplication by doing~$1-(1-x)\cdot(1+(1+x^2)\cdot(x+x^2))$.
This makes the evaluation of sequences of that have sizes that are power of five demand only~$4\cdot \operatorname{log}_5(N)-2 = 1.72\ldots\cdot \operatorname{log}_2(N)-2$, which is advantageous in comparison with the usual binary and ternary decomposition.
%
%
%
In Table~\ref{tab:algs} we provide smart factorizations for some prime sizes sequences with the complexity for the associated power of prime size sequences.
\begin{table*}[h]
\centering
\caption{Smart factorizations for evaluation of geometric series for prime size~$P$ and complexity for power of prime sizes~$N$.}
\begin{tabular}{ccc}\toprule
Size~$P$ & Algorithm for size~$P$& Complexity for size~$N$\\\midrule
\multirow{1}{*}{$P = 2$} & $f(2,x) = 1+x$ & \multirow{1}{*}{$2\cdot \operatorname{log}_2(N)-2$}\\\midrule
\multirow{1}{*}{$P = 3$} & $f(3,x) = 1+x+x^2$ & \multirow{1}{*}{$1.89\ldots\cdot \operatorname{log}_2(N)-2$}\\\midrule
\multirow{2}{*}{$P = 5$} & $y = x^2$ & \multirow{2}{*}{$1.72\ldots\cdot \operatorname{log}_2(N)-2$}\\
 & $f(5,x) = 1+(1+y)\cdot(x+y)$ & \\\midrule

\multirow{2}{*}{$P = 7$} & $y = x^2, w = y^2$ & \multirow{2}{*}{$1.78\ldots\cdot \operatorname{log}_2(N)-2$}\\
 & $f(7,x) = 1+(x+y)\cdot(1+y+w)$ & \\\midrule

\multirow{2}{*}{$P = 11$} & $y = x^2, w = y^2$ & \multirow{2}{*}{$1.73\ldots\cdot \operatorname{log}_2(N)-2$}\\
 & $f(11,x) = 1+(x+y)\cdot(1+(x+y)\cdot(1+w))$ & \\\midrule

\end{tabular}
\label{tab:algs}
\end{table*}
From Table~\ref{tab:algs}, the attentive reader might notice that these smart factorizations follow some underlying rule.
Indeed, that is what is happening.
For a prime size~$P$, we use the factorization
\begin{equation}
\scriptsize
\label{eq:mod2a}
f(P,x) = 
\begin{cases}
(1+x^2)\cdot f(\frac{P}{2},x^2),&\text{ if } P \equiv 0 \bmod 2 \\
1+(x+x^2)\cdot f(\frac{P-1}{2},x^2),&\text{ if } P \equiv 1 \bmod 2,
\end{cases}
\end{equation}
and we repeat the process interactively.
From~\eqref{eq:mod2a}, one can notice that it demand~$2\lfloor \operatorname{log}_2(P)\rfloor+t-2$, where~$t$ is the position of the second most significant bit in the binary representation of~$P$~\cite{Dimitrov1994,Lei1992}.

Also, note that once~$f(P,x)$ has been computed, we evaluate the next power of~$x$, namely~$x^P$, using the geometric series formula leading to~$x^P = 1+ f(P,x)\cdot(x-1)$.
%
%
This is recursively used in connection with Corollary~\ref{cor:reduction-geometric-series-prime}.

\begin{figure}[h]
\centering
\psfrag{PPPP/log}{$P'/\operatorname{log}_2(P)$}
\psfrag{P}{$P$}
\psfrag{1.7}{$1.7$}
\psfrag{1.9}{$1.9$}
\psfrag{2.1}{$2.1$}
\psfrag{2}{$2$}
\psfrag{6}{$6$}
\psfrag{10}{$10$}
\psfrag{14}{$14$}
\psfrag{18}{$18$}
\psfrag{21}{$21$}

\includegraphics[scale=0.9]{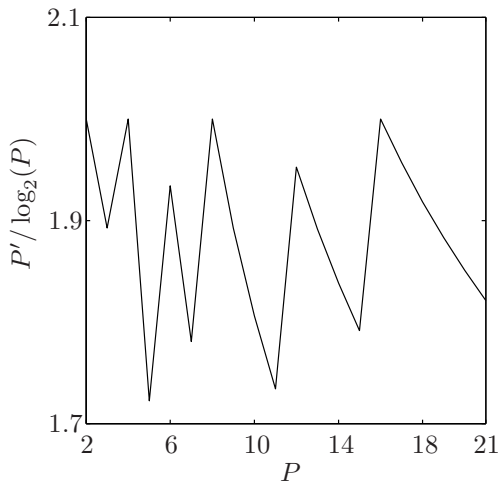}
\caption{Graphical representation of coefficient~$P'/\operatorname{log}_{2}(P)$ for the asymptotic complexity of sequences of size power of~$P$, where~$P' = 2\lfloor \operatorname{log}_2(P)\rfloor+t-2)$.}
\label{fig:graph-comp}
\end{figure}

%
In such a case, the computation of~$f(P,x)$ requires only~$P' < P$ multiplications, where~$P' = 2\lfloor \operatorname{log}_2(P)\rfloor+t-2$.
It implies that the result in Theorem~\ref{theo:complexity-direct-prime} changes to~$P'\cdot \operatorname{log}_P(N)-2 = P'/\operatorname{log}_2(P)\cdot \operatorname{log}_2(N)-2$.
%
%
Because of the importance of this result, we put in the format of Theorem~\ref{theo:complexity-fast-prime}.
\begin{theorem}[Complexity of Fast Evaluation of Power of Prime Size Series]
\label{theo:complexity-fast-prime}
Let~$N$ be a power of a prime~$P$.
%
%
Then, the computation of~$f(N, x)$ requires only~$P'\cdot \operatorname{log}_{P}(N)-2$ multiplications, where~$2\lfloor \operatorname{log}_2(P)\rfloor+t-2$ and~$t$ is the second most significant bit of~$P$ in binary representation.
\begin{proof}
Follow the proof for Theorem~\ref{theo:complexity-direct-prime} changing the number of multiplications required to compute~$f(P,x)$ from~$P$ to~$P' = 2\lfloor \operatorname{log}_2(P)\rfloor+t-2$.
\end{proof}
\end{theorem}
For a fair comparison, we should convert~$P'\cdot \operatorname{log}_{P}(N)-2$ to base~$2$, which result in~$P'/\operatorname{log}_{2}(P) \cdot \operatorname{log}_{P}(N)-2$.
These are the values shown on the third column of Table~\ref{tab:algs}.
Figure~\ref{fig:graph-comp} shows the coefficient~$P'/\operatorname{log}_{2}(P)$ which determine the asymptotic complexity in the case of sequences of size power of~$P$.
Note that base~$5$ apparently is a good candidate, as shown on Table~\ref{tab:algs} as well.
%
%

%
%

\section{Optimal Basis and Asymptotic Complexity}
\label{sec:optimal-base}
By the results shown on Table~\ref{tab:algs}, one may be tempted to conjecture that basis five is optimal.
Indeed, basis five is a good base, but it is not optimal.
One can see that by noting that~$f(25,x)$ can be computed by making~$f(25,x) = f(5,x)\cdot f(5,x^5)$.
Explicitly, we have
\begin{equation*}
f(25,x) = (1+x+x^2+x^3+x^4)\cdot(1+x^5+x^{10}+x^{15}+x^{20}).
\end{equation*}
Since the computation of both~$f(5,x)$ and~$f(5,x^5)$ requires only two multiplications, we have that~$f(25,x)$ requires six multiplications.
If one uses binary basis, it would be required~$2\operatorname{log}_2(25)-2 \approx 7$ multiplications.
So far it is not surprising that we saved one multiplication.
But have a look at the following computation of~$f(26,x)$:
\begin{equation*}
(1+x\cdot(1+x+x^2+x^3+x^4))\cdot(1+x^5+x^{10}+x^{15}+x^{20}).
\end{equation*}
It is equivalent to write~$f(26,x) = (1+x\cdot f(5,x))\cdot f(5,x^5)$.
One could say that it requires only one more multiplication than~$f(25,x)$, but it does not.
Indeed, if one follows 
\begin{align*}
y & = 1+(1+x^2)\cdot(x+x^2)\\
z & = 1+x\cdot y\\
v & = z-y\\
w & = 1+(1+v^2)\cdot(v+v^2)\\
t & = z\cdot w,
\end{align*}
it will need only six multiplications again due to the fact that~$f(26,x) = t$.
According Theorem~\ref{theo:complexity-fast-prime}, sequences that are power of~$26$ requires~$8\operatorname{log}_{26}(N)-2 = 1.7019 \cdots\operatorname{log}_{2}(N)-2$ multiplications, which is better than the complexity given by basis five and all the basis commented before.
Note that~$26 = 5^2+1$.
This is important and this pattern will follow for others sequences sizes.
One could use this approach with~$677 = 26^2+1$, leading to~$f(677,x) = (1+x\cdot f(26,x))\cdot f(26,x^{26})$.
For~$f(677,x)$, we would then need only fourteen multiplications:~six to compute~$f(26,x)$, six to compute~$f(26,x^{26})$ plus two multiplications.
For sequences that are power of~$677$, Theorem~\ref{theo:complexity-fast-prime} guarantees that we need only~$16\operatorname{log}_{677}(N)-2 = 1.7015 \cdots\operatorname{log}_{2}(N)-2$ multiplications, which is better than the complexity given by basis five, twenty-six and all the basis commented before.
However, note that the gain starts to be marginal and very little compared to basis twenty-six.
As one might have noticed, we are approaching the limit.
In the following, we show the asymptotic case, which must be of theoretical interest.

Consider the sequence~$y_{n} = y_{n-1}^2+1$, where~$y_0 = 1$.
Aho and Sloane~\cite{Aho1973} states the following theorem for recurrences.
\begin{theorem}
\label{theom:sloan}
Let~$y_n$ be defined by the quadratic recurrence sequence~$y_{n} = y_{n-1}^2+1$ with starting point~$y_0$.
We have the solution~$y_n = \lfloor k^{2^n} \rfloor$, where
\begin{equation*}
k = y_0 \exp\left( \sum_{n=0}^{\infty} 2^{-n-1} \operatorname{log}\left(1+y_n^{-2}\right) \right).
\end{equation*}
\end{theorem}
Note that this result is completely related to the following representations~$5_2 = (101), 26_5 = (101), 677 = (101)_{26}, \ldots$.
This is not a coincidence.
Indeed, it is a consequence of the series~$y_n = y_{n-1}^2+1$ with~$y_0 = 1$.
It is important because the series of length~$y_n$ can be evaluated with only twice the number of multiplications required by a series of size~$y_{n-1}$.
Theorem~\ref{theom:sloan} tells us that while the complexity is increasing exponentially, the size of the series is increasing~\emph{doubly} exponentially.
In other words, the size is growing faster than the complexity.

For the particular case treated here, note that the internal argument of the exponential function is a convergent series.
This can be seen using the geometric test~\cite{Graham1989}.
Being that said, using computer assistance we can estimate~$k = 1.50283680104976 \ldots$.

In order to estimate the complexity in the asymptotic case, note that the number of multiplication required by series that are power of the sizes~$\{1, 2, 5, 26, 677, \ldots\}$ requires~$\{0, 2\operatorname{log}_2(N)-2, 4\operatorname{log}_5(N)-2, 8\operatorname{log}_{26}(N)-2, 16\operatorname{log}_{677}(N)-2, \ldots\}$ multiplications, respectively.
In general, for series with sizes that are power of~$y_n$, we need to perform~$2^n\operatorname{log}_{y_n}(N)-2$.
Applying the result in Theorem~\ref{theom:sloan}, we have that
\begin{align*}
2^n\operatorname{log}_{y_n}(N)-2 & = 2^n\operatorname{log}_{k^{2^n}}(N)-2\\
& = \operatorname{log}_k(N)-2\\
& = \operatorname{log}_k(2)\operatorname{log}_2(N)-2\\
& = 1.70158214004473\ldots \operatorname{log}_2(N)-2.
\end{align*}
Note that this limit value is very close to the one obtained to the case of series with size power of~$677$, which is equals to the limit up to six decimal digits.
However, to achieve close to optimal performance in terms of complexity reduction, one must consider very large series, which may be impractical.

\section{Basis Combination}
\label{sec:base-combination}
Not always it is possible to decompose a geometric series as a product of prime size sequences of the same length.
This is particular true in cases where the original series size is not close to any number that is a power of a prime.
In order to reduce the multiplicative complexity and efficiently use the algorithms developed for small sizes, one can decompose the original series into the product of several others with different sizes using Theorem~\ref{theo:geometric-series-reduction}.

However, it may happen that the series size can be approximated to a series with size that is power of different integers, but not necessarily a power of a single prime number.
In this case, one could use directly~Theorem~\ref{theo:geometric-series-reduction}, or try to find a more efficient way to do that.
In this regard, let us consider the case where the series can be approximated with a number of terms that is  a product of powers different integers.
In order to do that, we need to consider~\emph{mixed basis}.
By mixed basis we mean a basis that is the combination of two or more prime integers.
Let us consider the starting example of basis compounded of~$\{2, 3\}$.
Using modular arithmetic representation for base~$2$ and~$3$, respectively, we have the cases outlined by~\eqref{eq:mod2} and~\eqref{eq:mod3}
\begin{figure*}
\begin{equation}
\label{eq:mod2}
f(N,x) = 
\begin{cases}
(1+x^2)\cdot f(\frac{N}{2},x^2),&\text{ if } N \equiv 0 \bmod 2 \\
1+(x+x^2)\cdot f(\frac{N-1}{2},x^2),&\text{ if } N \equiv 1 \bmod 2.
\end{cases}
\end{equation}
\end{figure*}
\begin{figure*}
\begin{equation}
\label{eq:mod3}
f(N,x) = 
\begin{cases}
(1+x+x^2)\cdot f(\frac{N}{3},x^3),&\text{ if } N \equiv 0 \bmod 3,\\
1+(x+x^2+x^3)\cdot f(\frac{N-1}{3},x^3),&\text{ if } N \equiv 1 \bmod 3,\\
1+x+(x^2+x^3+x^4)\cdot f(\frac{N-2}{3},x^3),&\text{ if } N \equiv 2 \bmod 3.
\end{cases}
\end{equation}
\end{figure*}
Let us have a closer look at the complexity for each case of base~$2$ and~$3$ decomposition.
For base~$2$, convert~$f(N,x)$ to either~$f(N/2,x^2)$ or~$f((N-1)/2,x^2)$ requires only~$2$ multiplications.
For base~$3$, if~$N \equiv 0 \bmod 3$, it is required only~$3$ multiplication to convert the evaluation of~$f(N,x)$ to~$f(N/3,x^3)$.
When~$N \equiv 1 \bmod 3$, we need~$3$ multiplications as well because of the decomposition~$1+(x+x^2+x^3)\cdot f((N-1)/3,x^3) = 1+x\cdot (1+x+x^2)\cdot f((N-1)/3,x^3)$. 
The computation of~$x^3$ can be done without multiplication because~$x^3 = x\cdot (1+x+x^2) - (x+x^2)$.
In these both case, where~$N \equiv 0, 1 \bmod 3$, use the base~$3$ is advantageous to using base~$2$, because~$3\operatorname{log}_3(N) < 2\operatorname{log}_2(N)$.
But this is not the case when~$N \equiv 2 \bmod 3$, because no matter we do, convert~$f(N,x)$ to~$1+x+(x^2+x^3+x^4)\cdot f((N-2)/3,x^3)$ requires at least~$4$ multiplications.
Let us estimate the average number of multiplications if we just use base~$3$.
For the sake of simplicity, let us suppose that we can pick a random number and following a uniform distribution~\cite{Ross2004}.
It is natural to assume that if we represent the picked number in ternary base, the chance of occurrence of the residues~$0$,~$1$ or~$2$ is also uniform.
In this case, the average number of multiplications for base~$3$ is~$(1/3 3+1/3 3+1/3 4)\operatorname{log}_3(N)-2 \approx 2.1\operatorname{log}_2(N)-2$.

The idea of combining base~$2$ with~$3$ is to use the decomposition induced by base~$3$ as much as possible, but adopt base~$2$ decomposition when~$N \equiv 2 \bmod 3$.
In this regard, we have the cases outlined by~\eqref{eq:mod6}, where~$r = \bmod 6$.
\begin{figure*}
\begin{equation}
\label{eq:mod6}
f(N,x) = 
\begin{cases}
(1+x+x^2)\cdot f(\frac{N-r}{3},x^3),&\text{ if } r = 0\text{ or }3,\\
1+(x+x^2+x^3)\cdot f(\frac{N-r}{3},x^3), & \text{ if } r = 1\text{ or }4,\\
(1+x^2)\cdot f(\frac{N-r}{2},x^2), & \text{ if } r = 2,\\
1+(x+x^2)\cdot f(\frac{N-r}{2},x^2), & \text{ if } r = 5.
\end{cases}
\end{equation}
\end{figure*}

\begin{figure}[h]
\centering
\includegraphics[scale=0.5]{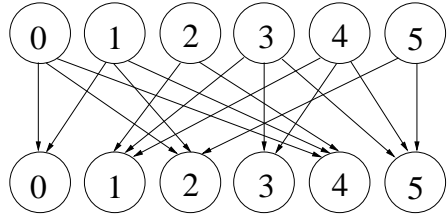}
\caption{Graphical representation of the transition between each iteration for the mixed basis combining~$2$ and~$3$.}
\label{fig:graphbase2and3}
\end{figure}

%
Supposing we pick again a random number according to uniform distribution, we want to evaluate the average number of multiplications required by a series with its size.
The evaluation is not straightforward as we have done for the previous case where there is no basis mixing.
Now we need to employ some stochastic analysis using Markov chains~\cite{Ross2004}.
For that, note that the way we decompose a series is a iterative process, and it depends on the value of~$N$.
Suppose we have first converted~$f(N,x)$ to~$f(N/3,x^3)$ because~$N \bmod 6 = 0$.
In the next iteration we will have to convert~$f(N/3,x^3)$ to something else, and the only possible values for~$N/3 \bmod 6$ are now~$0, 2$ or~$4$.
We build the Markov chain for this particular basis mixing based on the dependence of each iteration.
Analyzing each case, we come to the graph shown on Figure~\ref{fig:graphbase2and3}.
Note that the rule that governs is: if~$N$ has remaining~$0$ or~$1$ when divided by~$3$, we use base~$3$, if not, we use base~$2$.
%
%
%
%
%
%
%
%

Let us denote by~$p_j^{(i)}$ the probability of the Markov chain assume state~$i$ at iteration~$j$, where~$j \in \{0, 1, 2, 3, 4, 5\}$ are the remainders by six in each iteration.
For this particular case, we have the following system of transition equations to find the asymptotic solution of Markov chains:
\begin{align*}
\left[
\begin{rsmallmatrix}
p_0^{(i+1)}\\
p_1^{(i+1)}\\
p_2^{(i+1)}\\
p_3^{(i+1)}\\
p_4^{(i+1)}\\
p_5^{(i+1)}
\end{rsmallmatrix}
\right] & = 
\left[
\begin{rsmallmatrix}
\frac{1}{3} & \frac{1}{3} & 0 & 0 & 0 & 0\\
0 & 0 & \frac{1}{2} & \frac{1}{3} & \frac{1}{3} & 0\\
\frac{1}{3} & \frac{1}{3} & 0 & 0 & 0 & \frac{1}{2}\\
0 & 0 & 0  & \frac{1}{3} & \frac{1}{3} & 0\\
\frac{1}{3} & \frac{1}{3} & \frac{1}{2} & 0 & 0 & 0\\
0 & 0 & 0 & \frac{1}{3} & \frac{1}{3} & \frac{1}{2}
\end{rsmallmatrix}
\right]
\left[
\begin{rsmallmatrix}
p_0^{(i)}\\
p_1^{(i)}\\
p_2^{(i)}\\
p_3^{(i)}\\
p_4^{(i)}\\
p_5^{(i)}
\end{rsmallmatrix}
\right],
\end{align*}
with initial values
%
~$p_j^{(i)} = 1/6$ for all values of~$j$.
The solution for this particular case is
%
~$p_j^{(\infty)} = 1/10$ for~$j = 0, 3$ and~$p_j^{(\infty)} = 2/10$ for~$j = 1, 2, 4, 5$.
With these results we can estimate the multiplicative complexity for this mixed base.
The probability of using base~$3$ is~$p_0^{(\infty)}+p_1^{(\infty)}+p_3^{(\infty)}+p_4^{(\infty)} =6/10$, while the probability using base~$2$ is~$p_2^{(\infty)}+p_5^{(\infty)} = 4/10$.
Remember that when we use base~$3$ we use~$3$ multiplications, either the remainder is~$0$ or~$1$, and when we use base~$2$ we need~$2$ multiplications.
So, the average basis is~$b = 2^{\frac{4}{10}}\cdot 3^{\frac{6}{10}} \approx 2.55 \ldots$.
Thus, we need~$(3\cdot 6/10+2\cdot 4/10)\operatorname{log}_b(N)-2 = 1.9245\ldots \operatorname{log}_2(N)-2$.

This particular example shows the general procedure for mixing other basis--not demonstrated here.
One need to evaluate the number of multiplication required by each basis in each case and then select the appropriate basis to be used.
Then the analysis is done using the associated Markov chain resulting from the transition between iterations.
The average multiplicative complexity must take into account the asymptotic probabilities that are solution of the associated Markov chain.
The designer must have in mind that the basis mixing is not limited to only two different basis.
Any finite number of basis can be mixed, as long the purpose is to reduce the multiplicative complexity.
For completeness, we show the complexity of some particular basis mixing in Table~\ref{tab:base-mixing}.

\begin{table}[h]
\centering
\caption{Multiplicative complexity for computing~$f(N,x)$ using mixed basis}
\begin{tabular}{cc} \toprule
Mixed Basis & Multiplicative Complexity \\ \midrule 
 3,2  & $1.9245\ldots \operatorname{log}_2(N)$ \\ \midrule
 7,2  & $1.9057\ldots \operatorname{log}_2(N)$ \\\midrule 
 7,3,2  & $1.8749\ldots \operatorname{log}_2(N)$ \\\midrule 
 5,2  & $1.8554\ldots \operatorname{log}_2(N)$ \\ \midrule
 5,3,2  & $1.8299\ldots \operatorname{log}_2(N)$ \\ \midrule 
 7,5,3,2  & $1.8106\ldots \operatorname{log}_2(N)$ \\ \midrule 
 11,5,3,2  & $1.8036\ldots \operatorname{log}_2(N)$ \\\midrule 
 11,7,5,3,2  & $1.7932\ldots \operatorname{log}_2(N)$ \\ \bottomrule 
\end{tabular}
\label{tab:base-mixing}
\end{table}

\section{Application to Matrix Inversion Approximation}
\label{sec:app}
%
%
%
%

The advancement of digital computer design allowed the development of artificial image generation techniques, collectively called computer generated imagery (CGI)~\cite{Aubry2015}.
CGI has been employed in different areas such as flight simulation~\cite{Zeevi1990, Y1985}, entertainment~\cite{Hall1989} and architecture~\cite{Boyles2009}.
Several steps are needed to generate artificial images, starting at modelling to image rendering.
On the image rendering step, usually the computer has to simulate the interference of different light sources across a scene and be able to separate the different sources.
This process requires the inversion of matrices of arbitrary size that depend on the scene size and the resolution employing the CGI techniques.
These matrices are called light transport matrices.
Different from most of applications, it is not required only solve a linear system such as frequently occurs with general inverse problems~\cite{Ribes2008}.
%
%
Instead, it is required the actual matrix inversion due to its physical meaning and interpretation~\cite{Ng2011}.

Several papers have been considered a low complexity approach for finding the light transport matrix inverse~\cite{Chandraker2011, Ng2011, Bai2010, Lecocq2000}.
They share in common that their approach is based on Neumann series approximation~\cite{Stewart1998}.
Indeed, the community has resorted to approach the matrix inversion instead of computing the exact matrix inversion.
This is due to the fact that light transport matrices may reach sizes higher than~$10^{5}\times 10^{5}$~\cite{Ng2011}, which would be computationally intensive and demand too much time for execution.
In some scenarios, light transport matrices may be approximated by sparse matrices, which alleviate--but does not remove--the computational complexity.

%
%
%
%
%

%
%

Another particular application one could employ the optimization of Neumann series evaluation for inverting matrices is in massive multiple-input-multiple-output (MIMO) systems.
Several papers have considered the problem of inverting big matrices in this context~\cite{Shao2016, Wu2013, Tu2015, Prabhu2014}.
Most of the papers have considered the case where the number of antennas at the base station is considerably larger than the number of users. 
This lead to a context where small size series are enough for approximating the inverse matrix~$\mathbf{A}^{-1}$.
However, this may not be always the case, requiring a series with more terms to better approximate the inverse matrix~\cite{Prabhu2014}.
In this scenario we may employ optimized Neumann series evaluation, leading to lower complexity matrix inversion approximation.

The Neumann series is capable of approximating the inverse of a non-singular matrix by several matrix-matrix multiplications.
Let~$\mathbf{A}$ be a square matrix such that its eigenvalue magnitudes are confined to the interval~$(0,2)$.
One can approximate the inverse~$\mathbf{A}^{-1}$ using the series with~$N$ terms by the expression on~\eqref{eq:invAN}.
%
%
The number of terms~$N$ considered for approximation depend on the application and precision requirements.
On the limit, the approximation converge to the exact inverse matrix.
%
 
For this scenarios, one could consider the binary or ternary approach, which demands less multiplications than the usual Horner rule.
However, since in these applications small size series are used, such as~$5, 6, 7, 8$ and~$9$, we will employ the different factorizations that we have developed in Section~\ref{sec:fast-alg-small} and shown in Table~\ref{tab:algs} when adequate.
%
%
%
%
%
%
%
%
%
We generated matrices of sizes~$50, 100, 250$ and~$500$ for simulation, all possessing eigenvalues on the interval~$(0,2)$.
We reproduced~$1000$ replicates for each matrix size.
Table~\ref{tab:simulations} shows the time results for time and precision achieved for simulations.
\begin{table*}[h]
\begin{center}
\caption{Simulation results comparing fast algorithms for Neumann series sizes~$5, 6, 7, 8$ and~$9$ and direct evaluation in time (in seconds) with~$1000$ replicates}
\label{tab:simulations}
\scriptsize
\begin{tabular}{ccccccccc}\toprule
 \multirow{2}{*}{Series size} & \multicolumn{2}{c@{\quad}}{Matrix size $50$} & \multicolumn{2}{c@{\quad}}{Matrix size $100$} & \multicolumn{2}{c@{\quad}}{Matrix size $250$} &\multicolumn{2}{c@{\quad}}{Matrix size $500$} \\\cmidrule{2-9}
& Direct & Fast Alg. & Direct & Fast Alg. & Direct & Fast Alg. & Direct & Fast Alg.\\\midrule  
$N = 5$ & $1.75\cdot 10^{-4}$ &$8.41\cdot 10^{-5}$ &$7.28\cdot 10^{-4}$ &$3.32\cdot 10^{-4}$ &$6.11\cdot 10^{-3}$ &$2.87\cdot 10^{-3}$ &$5.5\cdot 10^{-2}$ &$2.16\cdot 10^{-2}$ \\\midrule
$N = 6$ & $2.04\cdot 10^{-4}$ &$1.07\cdot 10^{-4}$ &$8.62\cdot 10^{-4}$ &$4.53\cdot 10^{-4}$ &$7.51\cdot 10^{-3}$ &$3.92\cdot 10^{-3}$ &$6.59\cdot 10^{-2}$ &$3.24\cdot 10^{-2}$ \\\midrule
$N = 7$ & $2.36\cdot 10^{-4}$ &$1.1\cdot 10^{-4}$ &$1.01\cdot 10^{-3}$ &$4.66\cdot 10^{-4}$ &$9.91\cdot 10^{-3}$ &$4.72\cdot 10^{-3}$ &$8.11\cdot 10^{-2}$ &$3.61\cdot 10^{-2}$ \\\midrule
$N = 8$ & $2.65\cdot 10^{-4}$ &$1.36\cdot 10^{-4}$ &$8.91\cdot 10^{-4}$ &$4.47\cdot 10^{-4}$ &$1.08\cdot 10^{-2}$ &$5.73\cdot 10^{-3}$ &$8.78\cdot 10^{-2}$ &$4.37\cdot 10^{-2}$ \\\midrule
$N = 9$ & $2.96\cdot 10^{-4}$ &$1.33\cdot 10^{-4}$ &$8.61\cdot 10^{-4}$ &$3.77\cdot 10^{-4}$ &$1.26\cdot 10^{-2}$ &$6.07\cdot 10^{-3}$ &$1.02\cdot 10^{-1}$ &$4.67\cdot 10^{-2}$ \\\bottomrule
\end{tabular}
\end{center}
\end{table*}
Note that the average time spent for computing the approximate inverse matrix using the algorithms on sections~\ref{sec:decomp-series-prime} and~\ref{sec:base-combination} are up to~$2.2$ times smaller when compared to the usual approach.
%
%
We employed usual binary and ternary factorizations for the approximations with series of sizes~$8$ and~$9$, respectively.
For the remaining sizes, we employed our proposed factorizations.

\section{Conclusions}
\label{sec:conclusion}
In this paper we have derived different factorizations for the evaluation of Neumann series.
The factorizations achieved reduced complexity--reducing the well-known binary decomposition with~$2\operatorname{log}_2(N)-2$ to the asymptotic case~$1.7\ldots \operatorname{log}_2(N)-2$.
The proposed factorization could achieve reduced complexity not even for asymptotic cases, but for small and average size series, such as powers of~$5, 7, 11$ and several others.
We have shown also that is possible to mix different basis and achieve further complexity reduction while taking advantage offered by each employed basis.
The complexity analysis was performed with the help of Markov chain theory after a examination of case by case of each state on the mixed basis factorization.
The theory exposed here was applied for the problem of approximate matrix inversion in the context of light transport matrix and massive MIMO systems.
Simulations comparing the proposed algorithms with a direct evaluation of Neumann series was performed.
We show that the proposed algorithms can drastically reduced the overall inversion time.
For small size matrix we can perform a inverse matrix approximation up to~$2.2$ times faster than the usual approach (series with size~$9$) in scenario where binary approach could not be applied.

Future works may include the asymptotic complexity in deeper detail.
Starting from a small length other than~$5$ that presents smaller complexity result and that leads to a smart factorization using the sequence~$y_n = y_{n-1}^2+1$ with~$y_0 = 1$ could lead to new lower bounds for asymptotic complexity.
Also, a possible different approach for analyzing the multiplicative complexity for mixed basis other than using Markov chains could be pursued.
%
%


\end{document}